\documentclass[10pt]{amsart}

\usepackage{amsmath, amsthm, amscd, amssymb, amsfonts, pstricks}

\input amssym.def
\input amssym.tex
\newtheorem{thm}{Theorem}[section]
\newtheorem*{thm*}{Theorem}

\newtheorem*{cor*}{Corollary}
\newtheorem{lem}[thm]{Lemma}
\newtheorem{prop}[thm]{Proposition}

\newtheorem*{con*}{Conjecture}
\newtheorem*{prob*}{Problem}

\theoremstyle{definition}

\theoremstyle{remark}
\newtheorem{rem}[thm]{Remark}

\newcommand{\hol}{{\text{Hol}}}
\newcommand{\aut}{{\text{Aut}}}

\newcommand{\cF}{{\mathcal{F}}}
\newcommand{\cC}{\mathcal{C}}

\newcommand{\cB}{\mathcal{B}}

\newcommand{\fin}{\mathcal{F}in}
\newcommand{\fbc}{\mathcal{F}\mathcal{B}\mathcal{C}}
\newcommand{\vc}{\mathcal{V}\mathcal{C}}
\newcommand{\all}{\mathcal{A}\mathcal{L}\mathcal{L}}

\def\<{\langle}
\def\>{\rangle}
\def\hol{{\rm Hol}}
\def\aut{{\rm Aut}}
\def\inn{{\rm Inn}}
\def\gl{{\rm GL}_2({\mathbb Z})}

\def\Z{\mathbb Z}
\def\K{\mathbb K}
\def\ar{\rightarrow}
\def\t{\tau}
\def\x{\xi}
\def\z{\zeta}

\begin{document}

\title[On the vanishing of the lower K-theory of $\hol(F_2)$]
{On the vanishing of the lower K-theory of the holomorph of a free group on two generators}

\author{Vassilis Metaftsis}	
\address{University of the Aegean
School of Sciences
Department of Mathematics
Gr-83200 Karlovassi, Samos, Greece}
\email{vmet@aegean.gr}

\author{Stratos Prassidis}
\address{Department of Mathematics and Statistics, Canisius College,
Buffalo, New York 14208, U.S.A}
\email{prasside@canisius.edu}

\begin{abstract}
We show that the holomorph of the free group on two generators satisfies the Farrell--Jones
Fibered Isomorphism Conjecture. As a consequence, we show that the lower $K$-theory of the above group vanishes.
\end{abstract}
\maketitle

\section{Introduction}

Certain obstructions that appear in problems of topological rigidity of manifolds are elements of algebraic $K$-groups, specially lower $K$-groups. For this reason, the calculation of the lower $K$-groups has implications in geometric topology.

The main modern tool for calculating lower $K$-groups (and other geometrically interesting obstruction groups) is the Farrell--Jones Fibered Isomorphism Conjecture. The Conjecture provides an inductive method for calculating the obstruction groups of a group from those of certain subgroups. More specifically, if a group satisfies the Fibered Isomorphism Conjecture for a specific theory, then  the obstruction groups can be calculated from the obstruction groups  of the {\it virtually cyclic subgroups}. The last class of subgroups consists the finite subgroups and groups that are virtually infinite cyclic. 
The virtually infinite cyclic groups are of two types:
\begin{itemize}
\item Groups $V$ that surject onto the infinite cyclic group $\Z$ with finite kernel, i.e. $V = H \rtimes \Z$ with $H$ finite.
\item Groups $W$ that surject onto the infinite dihedral group $D_{\infty}$ with finite kernel, i.e.
$W  = A*_BC$ with $B$ finite and $[A:B] = [C:B] = 2$.
\end{itemize}
Obstruction groups that can be calculated this way are pseudoisotopy groups, $K$-groups of group rings, $K$-groups of $C^*$-algebras, $L^{-\infty}$-groups.

The fundamental work of Farrell--Jones (\cite{fj1}) deals with the Fibered Isomorphism Conjecture for the pseudoisotopy spectrum. It should be noted that if a group satisfies the Fibered Isomorphism Conjecture for pseudoisotopies, then it satisfies the Isomorphism Conjecture for the lower $K$-groups. The reason for this is that the lower homotopy groups of the pseudoisotopy spectrum and the $K$-theory spectrum are isomorphic. For $K$-groups there is a refinement of the Conjecture that was given in \cite{dkr}.
They showed that finite and  virtually infinite cyclic subgroups of the first type suffice in detecting the $K$-theory of the group.

Our main interest is in computing the lower $K$-groups of the holomorph of $F_2$, the free group on two generators. For a group $G$, the holomorph of $G$ is the universal split extension of $G$. Thus, it fits into a split exact sequence:
$$1 \to G \to \hol(G) \to \aut(G) \to 1.$$ 
The main result of the paper is the following:

\begin{thm*}[Main Theorem]
The group $\hol(F_2)$ satisfies the Fibered Isomorphism Conjecture for pseudoisotopies. Furthermore, if $\Gamma < \hol(F_2)$ then 
$$\text{Wh}({\Gamma}) = \widetilde{K}_0({\Z}{\Gamma}) = K_i({\Z}{\Gamma}) = 0, \quad\text{for}\;\;
i \le -1.$$
\end{thm*}

Notice that $\hol(F_2)$ is equipped with a sequence of surjections:
$$\hol(F_2) \to \aut(F_2) \to \gl.$$
The second surjection is induced by sending an automorphism of $F_2$ to an automorphism of its abelianization ${\Z}^2$. Notice that the kernel of both surjections is isomorphic to $F_2$.
To show that $\hol(F_2)$ satisfies the Fibered Isomorphism Conjecture, we use the fact that every automorphism of $F_2$ is {\it geometric}, i.e. it can be realized by a diffeomorphism on a surface with boundary. That allows us to show first that $\aut(F_2)$ satisfies the Conjecture and using the same fact again that $\hol(F_2)$ does also.

The group $\gl$ splits as an amalgamated free product of finite dihedral groups:
$$\gl = D_4*_{D_2}D_6.$$
Thus both $\aut(F_2)$ and $\hol(F_2)$ split as amalgamated free products.
Using this fact and the properties of elements of finite order in $\aut(F_2)$ given in \cite{meskin,dd}, we determine the finite and the virtually infinite cyclic subgroups of $\aut(F_2)$ and $\hol(F_2)$, up to isomorphism. The list of groups is short and their lower $K$-theory vanishes. The Main Theorem follows from this observation.

\bigskip

{\it The authors would like to thank Tom Farrell for asking the question on the lower $K$-theory of $\aut(F_2)$.}

\section{Preliminaries}

Let $G$ be a discrete group. By a {\it class of subgroups} of $G$ we wean a collection of subgroups of $G$ that is closed under taking subgroups and conjugates. In our application we consider the following classes of subgroups:
\begin{itemize}
\item ${\cF}in$, the class of finite subgroups of $G$.
\item ${\cF}{\cB}{\cC}$ for the class of finite by cyclic subgroups. Those are subgroups $H < G$ such that
$$1 \to A \to H \to C \to 1$$
where $C$ is cyclic (finite or infinite) group and $A$ is finite. Notice that ${\cF}{\cB}{\cC}$ contains 
${\cF}in$ and subgroups $H = A {\rtimes}{\Z}$, when $C$ is the infinite cyclic group.
\item $\vc$ for the class of virtually cyclic subgroups of $G$. For a such a subgroup $H$ either 
$H \in \fin$, $H \in \fbc$, or $H = A*_BC$, where $A$, $B$, $C$ are finite and $[A:B] = [C:B] = 2$.
\item $\all$ for the class of all subgroups of $G$.
\end{itemize}

Let $\cC$ be a class of subgroups of $G$. The {\it classifying space} for $\cC$, $E_{\cC}G$ is a $G$-CW-complex such that the isotropy groups of the actions are in $\cC$ and, for each $H \in \cC$, the fixed point set of $H$ is contractible (for more details \cite{dl}, \cite{luck}). 

\begin{rem}
For $\aut(F_2)$ the classifying space for finite groups is the auter space (\cite{hv}, \cite{jensen}).
\end{rem}

The Fibered Isomorphism Conjecture (FIC) was stated by Farrell--Jones (\cite{fj1}). For the groups that holds, it provides an inductive method for computing obstruction groups in geometric topology (for a review see \cite{lr}). If $G$ satisfies the FIC then the natural map
$$H_n^G(E_{\vc}G; \K\Z^{-\infty}) \to H_n^G(E_{\all}G; \K\Z^{-\infty}) = K_n(\Z G)$$
is an isomorphism for $n \le 1$. Notice that the left hand side of the isomorphism can be computed from the virtually cyclic subgroups of $G$.

In general, there are ``forgetful maps''
$$H_n^G(E_{\fin}G; \K\Z^{-\infty}) \to H_n^G(E_{\fbc}G; \K\Z^{-\infty}) \to
H_n^G(E_{\vc}G; \K\Z^{-\infty})$$
The difference between the class $\fin$ and the class $\vc$ is that the second class can be captured by  the Waldhausen and Bass--Farrell Nil-groups of the infinite virtually cyclic subgroups.
In \cite{davis} and \cite{dkr}, it was shown that the second map is an isomorphism. Essentially, the authors proved that the Waldhausen's Nil-groups that appear in the $K$-theory of virtually infinite cyclic subgroups can be detected by the Bass--Farrell Nil-groups that appear in the $\fbc$ class. 

The FIC is known to hold for certain classes of groups. One class of interest for this paper is the class of {\it strongly poly-free groups}. A group $\Gamma$ is called strongly poly-free if there is a filtration:
$${\Gamma} = {\Gamma}_0 \ge {\Gamma}_1 \ge \dots \ge {\Gamma}_n = \{1\}.$$
such that:
\begin{enumerate}
\item  ${\Gamma}_i$ is normal in $\Gamma$ for each $i$.
\item ${\Gamma}_{i}/{\Gamma}_{i+1}$ is finitely generated free for all $0 \le i \le n - 1$.
\item For each ${\gamma}\in \Gamma$ there is a compact surface $S$ and a diffeomorphism 
$f: S \to S$ such that the induced homomorphism $f_*$ on 
${\pi}_1(S)$ is equal to $c_{\gamma}$ in $\text{Out}({\pi}_1(S))$, where $c_{\gamma}$ is the action of 
$\gamma$ 
on  ${\Gamma}_{i}/{\Gamma}_{i+1}$ by conjugation and ${\pi}_1(S)$ is identified with 
${\Gamma}_{i}/{\Gamma}_{i+1}$ via a suitable isomorphism.
\end{enumerate}
In \cite{afr} and \cite{fr} it was shown that a finite extension of a strongly poly-free group satisfies the FIC.

\begin{rem}\label{rem-strongly}
Let $\Gamma$ be a group that satisfies (1) and (2) above. We assume that ${\Gamma}_{i}/{\Gamma}_{i+1} \cong F_2$. Then $G$ is strongly poly-free. For this, let $T^2$ be the torus and $p = (1, 1)\in T^2$. Then 
${\pi}_1(T^2 \setminus \{p\}, x) = F_2$. In this case, 
$$\text{Out}(F_2) = \aut(F_2)/\inn(F_2) = \gl$$
where $\aut(F_2)$ denotes the automorphism group of $F_2$. 
Let $c_{\gamma}$ be an induced homomorphism as in (3) above. Then the image of $c_{\gamma}$ to 
$\text{Out}(F_2)$ can be represented by a diffeomorphism $f$ of $T^2$ that fixes $p$. After an isotopy starting at the identity on $T^2$, we can assume that $f$ fixes a small open disk $D$ around $p$. 
Then $f$ induces a diffeomorphism on the compact surface
$$f: T^2\setminus D \to T^2 \setminus D,$$
that fixes the boundary. Thus $f_* = c_{\gamma}$ in $\text{Out}(F_2)$.
\end{rem}

Start with an exact sequence of groups.
$$1 \to A \to B \xrightarrow{r} C \to 1.$$
In the Appendix of \cite{fj1}, it was shown the FIC holds for $B$ if:
\begin{itemize}
\item It holds for $C$.
\item For each virtually cyclic subgroup $V$ of $C$, it holds for $r^{-1}(V)$.
\end{itemize}

Using this result we show the following

\begin{prop}\label{prop-f2}
Let
$$1 \to F_2 \to G \xrightarrow{r} H \to 1$$
be an exact sequence. If the FIC holds for $H$, then it holds for $G$.
\end{prop}

\begin{proof}
Using the result in \cite{fj1}, it is enough to show that the FIC holds for $r^{-1}(V)$, where $V$ is a virtually cyclic subgroup of $H$. 

If $V$ is finite, then $r^{-1}(V)$ is a finite extension of $F_2$ and $r^{-1}(V)$ is  a finite extension 
of a free group. The result follows from Remark \ref{rem-strongly}.

If $V$ is infinite, then $V$ contains an infinite cyclic normal subgroup $W$ of finite index. Then 
$r^{-1}(W)$ is a normal subgroup of $r^{-1}(V)$ and fits into an exact sequence:
$$1 \to F_2 \to r^{-1}(W) \to W \to 1.$$
Then there is a filtration $r^{-1}(W) > F_2 > \{1\}$, with the first quotient being an infinite cyclic group. Obviously, every homomorphism of $\Z$ is realized by a diffeomorphism of $S^1{\times}[0, 1]$. 
Using Remark \ref{rem-strongly}, we see that $r^{-1}(W)$ is strongly poly-free. Therefore, $r^{-1}(V)$ is a finite extension of a strongly poly-free group. By \cite{fr}, it satisfies the FIC, completing the proof of the proposition.
\end{proof}

Let $\hol(F_2)$ denote the holomorph of $F_2$, namely, the universal split extension of $F_2$:
$$1 \to F_2 \to \hol(F_2) \xrightarrow{p} \aut(F_2) \to 1.$$
Notice that there is an exact sequence
$$1 \to \text{Inn}(F_2) \to \aut(F_2) \xrightarrow{q} \gl \to 1$$
that is induced by mapping the automorphisms of $F_2$ to the automorphisms of its abelianization.
That induces an exact sequence:
$$1 \to F_2 \to \aut(F_2) \xrightarrow{q} \gl \to 1.$$

\begin{prop}\label{prop-hol}
The FIC holds for $\aut(F_2)$ and $\hol(F_2)$.
\end{prop}

\begin{proof}
The group $\gl$ contains a subgroup of finite index that is isomorphic to $F_2$. In fact,
the following short exact sequence is known to hold, as a result of the standard action that
$\gl$ admits on the upper half plane:
$$1\to F_2\to \gl\to D_{12}\to 1$$
(see for example \cite{dd}). 
Thus the FIC holds for $\gl$. Now $\aut(F_2)$ fits into an exact sequence:
$$1 \to F_2 \to \aut(F_2) \to \gl \to 1.$$
By Proposition \ref{prop-f2}, the FIC holds for $\aut(F_2)$. Also, $\hol(F_2)$ fits into an exact sequence:
$$1 \to F_2 \to \hol(F_2) \to \aut(F_2) \to 1.$$
By Proposition \ref{prop-f2} again, the FIC holds for $\hol(F_2)$.
\end{proof}

\section{Infinite Finite-by-Cyclic Subgroups of $\hol(F_2)$}

Since there is an exact sequence
$$1 \to \text{Inn}(F_2) \to \aut(F_2) \xrightarrow{p} \gl \to 1$$
and $F_2=\< a,b\>$ is torsion free, every finite subgroup of $\aut(F_2)$ maps
isomorphically to a finite subgroup of $\gl$. On the other hand,
$\gl$ admits a decomposition as an amalgamated free product of the form 
\begin{equation}\label{decompgl}
\gl = D_4 *_{D_2} D_6
\end{equation}
where $D_2$, $D_4$ and $D_6$ are dihedral groups of orders $4$, $8$ and $12$
respectively. Hence, any finite subgroup of $\gl$ is a subgroup of a conjugate
of either $D_2$ or $D_4$ or $D_6$ and hence, so is every finite subgroup
of $\aut(F_2)$.

Now a presentation for $\aut(F_2)$ is given by
$$\< p, x,y,\t_a,\t_b\mid x^4=p^2=(px)^2=1, (py)^2=\t_b, x^2=y^3\t_b^{-1}\t_a,$$
$$p^{-1}\t_ap=x^{-1}\t_ax=y^{-1}\t_ay=\t_b, p^{-1}\t_bp=\t_a, x^{-1}\t_bx=\t_a^{-1}, y^{-1}\t_by=\t_a^{-1}\t_b\>$$
where $\t_a,\t_b$ are the inner automorphism of $F_2$ corresponding to $a,b$
respectively (see for example \cite{mks}). Moreover, a presentation for $\gl$ is
given by
$$\gl=\< P,X,Y\mid X^4=P^2=(PX)^2=(PY)^2=1, X^2=Y^3\>$$ and $G$ maps
onto $\gl$ by $p\mapsto P$, $x\mapsto X$, $y\mapsto Y$, $\t_a,\t_b\mapsto 1$.

As shown in \cite{meskin}, if $g$ is an element of finite order
in $\aut(F_2)$, then $g$ is conjugate in $\aut(F_2)$, to one of the following elements
$p, px, px\t_a, x^2, y^2\t_b^{-1}$ or $x$ with orders $2,2,2,2,3$ or $4$
respectively. This fact implies that $\aut(F_2)$ cannot contain finite
subgroups isomorphic to $D_6$.   
Moreover, any element of $\aut(F_2)$ can be written uniquely
in the form $p^{r}u(x,y)x^{2s}w(\t_a,\t_b)$ where $r,s\in \{0,1\}$, $w(\t_a,\t_b)$ is a reduced
word in $\inn(F_2)$ and $u(x,y)$ is a reduced word where $x,y,y^{-1}$ are the only
powers of $x,y$ appearing (see \cite{meskin,mks}). 

Also, due to the decomposition (\ref{decompgl}) of $\gl$, $\aut(F_2)$ is
also an amalgamated free product of the form 
\begin{equation}\label{decompaut}
\aut(F_2)=B*_{D}C
\end{equation}
where $B,C$ and $D$ fit into the following short exact sequences
$$1\to \inn(F_2)\to B\to D_4\to 1$$
$$1\to \inn(F_2)\to C\to D_6\to 1$$
$$1\to \inn(F_2)\to D\to D_2\to 1.$$

Moreover, since every one of $B,C$ and $D$ are free-by-finite groups they admit
an action on a tree with finite quotient graph and finite vertex and edge
stabilizers (as a corollary of the Almost Stability Theorem of Dicks and 
Dunwoody \cite{dd}). In fact, they are also amalgamated free products of the
form 
\begin{eqnarray}\label{furtherdecomp}
B &=& D_4*_{\Z/2\Z}D_2 = \< x,p\> *_{\< px\>} \< px, x^2\t_b\>\nonumber\\
C &=& D_3*_{\Z/2\Z}D_2 = \<y^2\t_b^{-1}\t_a, p\> *_{\< p\>} \< p,x^2\>\\
D &=& D_2*(\Z/2\Z) = \< p,x^2\> * \< x^2\t_b\>.\nonumber
\end{eqnarray}
Once again, the elements of finite order are $p,px, x^2\t_b, px^3\t_b, x^2, px^2, y^2\t_b^{-1}\t_a$ and $x$.
To be in accordance with Meskin, we see that $x^2 (y^2t_b^{-1}\t_a) x^2=y^3\t_b^{-1}\t_a (y^2\t_b^{-1}\t_a)
\t_a^{-1}\t_by^{-3}=y^2\t_b^{-1}$, 
$\t_a^{-1}x^{-1}(px^3\t_b)x\t_a=px\t_a$, $y(x^2\t_b)y^{-1}=x^2$ and that $x^{-1}(px^2)x=p$.  
 
By definition, $G=\hol(F_2)$ is the universal 
split extension of $\aut(F_2)$ and thus it fits to the split exact sequence
$$1\ar F_2 \ar \hol(F_2)\ar \aut(F_2)\ar 1.$$
So $\hol(F_2) = F_2 \rtimes \aut(F_2)$.  Hence, the above presentation
for $\aut(F_2)$ provides us with a presentation for $\hol(F_2)$. Namely,
$$\hol(F_2)=\< p,x,y,\t_a,\t_b,a,b\mid x^4=p^2=(px)^2=1, (py)^2=\t_b, x^2=y^3\t_b^{-1}\t_a,$$
$$p^{-1}\t_ap=x^{-1}\t_ax=y^{-1}\t_ay=\t_b, p^{-1}\t_bp=\t_a, x^{-1}\t_bx=\t_a^{-1}, y^{-1}\t_by=\t_a^{-1}\t_b,$$
$$\t_a^{-1} a\t_a=a, \t_a^{-1} b \t_a=a^{-1}ba, \t_b^{-1}a\t_b=b^{-1}ab, \t_b^{-1}b\t_b=b,
p^{-1}ap=b, p^{-1}bp=a,$$
$$x^{-1}ax=b, x^{-1}bx=a^{-1}, y^{-1}ay=b, y^{-1}by=a^{-1}b\>.$$

Moreover, the decomposition (\ref{decompaut}) of $\aut(F_2)$ provides an amalgamated
free product decomposition for $\hol(F_2)$:
\begin{equation}\label{simplehol}
\hol(F_2)= (F_2\rtimes B) *_{F_2\rtimes D} (F_2\rtimes C)
\end{equation}
and based on (\ref{furtherdecomp}) we have 
\begin{eqnarray}\label{furtherhol}
F_2\rtimes B  &=& (F_2\rtimes D_4)*_{F_2\rtimes \Z/2\Z}(F_2 \rtimes D_2)\nonumber\\
F_2\rtimes C &=&  (F_2\rtimes D_3)*_{F_2\rtimes \Z/2\Z}(F_2 \rtimes D_2)\\
F_2\rtimes D &=& (F_2\rtimes D_2)*(F_2\rtimes \Z/2\Z)\nonumber
\end{eqnarray}

Based again on the Almost Stability Theorem, we see that every vertex group in the 
above graphs of groups is a free-by-finite group so, it also admits a decomposition 
as a graph of groups with finite vertex groups. An analysis, based on the
presentations and also on the fact that the action of $x,y,p$ on $a,b$ is the same as that on $\t_a,\t_b$, would give us the following:
In $F_2\rtimes B$, 
$$F_2\rtimes D_4=D_4 *_{\Z/2\Z} D_2 = \< x,p\>*_{\< px\>}\< px, x^2b\>$$
$$F_2\rtimes D_2=D_2*_{\Z/2\Z}D_2*_{\Z/2\Z}D_2 = \< px ,x^2\t_b\>*_{\<px\>}\<px, x^2\t_bb^{-1}\> *_{\< x^2\t_bb^{-1}\>} \< pxa, x^2\t_bb^{-1}\>$$
$$F_2 \rtimes \Z/2\Z = (\Z/2\Z*_{b}) *\Z/2\Z = (\< px\>*_{b}) * \< pxa\>.$$    
In $F_2\rtimes C$,
$$F_2\rtimes D_3= D_3 *_{\Z/2\Z} D_3 = \< y^2\t_b^{-1}, p\> *_{\< p \>} 
\< y^2\t_b^{-1}a, p \>$$
$$F_2\rtimes D_2= D_2 * \Z/2\Z= \< p,x^2 \> * \<x^2b\>$$
$$F_2 \rtimes \Z/2\Z = (\Z/2\Z *\Z/2\Z)*_{a} = (\<p\> * \<pba^{-1}\>)*_{a}.$$
Finally, in $F_2\rtimes D$,
$$F_2\rtimes D_2= D_2 * \Z/2\Z = \< p, x^2\> * \< x^2b\>$$
$$F_2\rtimes \Z/2\Z = \Z/2\Z * \Z/2\Z *\Z/2\Z = \<x^2\t_b\> * \< x^2\t_bb\> *\< x^2\t_bb^{-1}a\>.$$
In the above, $S*_t$ denotes the HNN-extension with base group $S$ and stable letter $t$.

For example, $F_2\rtimes \< px, x^2\t_b\>$ has a presentation of the form
$$\< \x_1, \x_2, a, b\mid \x_1^2=\x_2^2=1, [\x_1,\x_2]=1, \x_1a\x_1=a^{-1}, \x_1b\x_1=b,
\x_2 a\x_2=b^{-1}a^{-1}b, \x_2b\x_2=b^{-1}\>$$ where $\x_1=px$ and $\x_2=x^2\t_b$.
By setting $\z_2=b\x_2$ and eliminating $b$, we get
$$\< \x_1,\x_2,a,\z_2\mid \x_1^2=\x_2^2=\z_2^2=1, [\x_1,\x_2]=1, \x_1a\x_1=a^{-1}, \z_2a\z_2=a^{-1},
[\x_1,\z_2]=1\>.$$ 
Now by setting $\x_3=\x_1a$ and eliminating $a$ we get
$$\< \x_1,\x_2,\x_3,\z_2 \mid \x_1^2=\x_2^2=\z_2^2=\x_3^2=1, [\x_1,\x_2]=[\z_2,\x_3]=[\x_1,\z_2]=1\>$$
which is the desired decomposition. 

So now we can prove the following result which generalizes the result in \cite{meskin} on the elements of finite order in $\aut(F_2)$.

\begin{lem}\label{lem-finite}
An element of finite order in $\hol(F_2)$ is conjugate to exactly one of $p,px,pxa, px\t_a, px\t_aa,
x^2, x^2b, y^2\t_b^{-1}, y^2\t_b^{-1}a$ and $x$ with orders $2$, $2$, $2$, $2$, $2$, $2$, $2$, $3$, $3$ and $4$ respectively.
\end{lem}

\begin{proof}
Given the above decomposition, every element of finite order is a conjugate of an element
of a vertex group. So it suffices to observe the following: $x^2 (px^3\t_bb) x^2= px\t_aa$, 
$x^2 (px^3\t_b) x^2= px\t_a$, $x^2(px^3b)x^2=pxa$, $x (px^2) x^{-1}=p$, $y(x^2\t_b)y^{-1}=x^2$,
$b^{-1}(x^2\t_bb^{-1})b=x^2\t_bb$, $x\t_a^{-1}\t_by^{-1}
(x^2\t_bb^{-1}a) y\t_b^{-1}\t_ax^{-1}=x^2\t_bb$ and $yb^{-1}x\t_a^{-1}\t_by^{-1}(x^2\t_bb)y\t_b^{-1}\t_ax^{-1}by^{-1}=x^2b$. Notice also that
$x^2b$ is no longer conjugate to $x^2$ since the relation $x^2=y^3\t_b^{-1}\t_a$ has no
equivalent for $a$ and $b$ due to the semidirect product structure of $G$. 
\end{proof}

From the fact that $\hol(F_2)=\< a,b\> \rtimes \aut(F_2)$ we have that every element $W$ of
$\hol(F_2)$ can be written uniquely in the form 
$$W=Vz(a,b)$$ 
where $V\in\aut(F_2)$ and $z(a,b)$ is
a word in the free group $\< a,b\>$. So, the normal form for the elements of
$\aut(F_2)$ implies the existence of a normal form for the elements of $\hol(F_2)$:
$$W=p^ru(x,y)x^{2s}w(\t_a,\t_b)z(a,b)$$ 
where $w(t_a,t_b)$ is a reduced word in the free
group $\< t_a, t_b\>$, $u(x,y)$ is a reduced word where $x,y,y^{-1}$ are the only
powers of $x,y$ appearing and $r,s\in\{ 0,1\}$. Moreover, every vertex group in
the decomposition (\ref{simplehol}) has also a normal form. More specifically,
every element in $F_2\rtimes B$ can be written uniquely in the form
$p^rx^nx^{2s}w(\t_a,\t_b)z(a,b)$ where $r,n,s\in \{0,1\}$ and every element
in $F_2\rtimes C$ can be written uniquely in the form
$p^ry^nx^{2s}w(\t_a,\t_b)z(a,b)$ where $r,s\in \{0,1\}$ and $n\in\{ 0,1,-1\}$
and $w(\t_a,\t_b)$ is a reduced word in $\<\t_a,\t_b\>$ and $z(a,b)$
is a reduced word in $\< a,b\>$.

Notice ]that one can define a natural epimorphism $$\hol(F_2) \ar \gl$$ with
kernel $\<a,b\>\rtimes \< \t_a,\t_b\>$. In fact, $\hol(F_2)$ fits into the following
short exact sequence
$$1\ar F_2\rtimes F_2\ar \hol(F_2) \ar \gl \ar 1$$
although such a sequence does not split.

We are searching for subgroups of $G=\hol(F_2)$ which are isomorphic
to $ A\rtimes \Z$ where $A$ is a finite subgroup of $G$. In our
argument we shall make extensive use of the following well known
result of Bass-Serre theory \cite{serre}. Let $M$ be a group that acts
on its standard tree $T$ and $m\in M$ such that $m$ stabilizes two distinct
vertices of $T$. Then $m$ stabilizes the (unique reduced) path that connects
the two vertices. In particular, $m$ is an element of every edge stabilizer
of every edge that constitute the path that connects the two vertices.
  
In fact we shall show the following:

\begin{prop}\label{prop-fbc}
The only finite-by-cyclic subgroups of $G$ are  ${\Z}/2{\Z}{\times}{\Z}$.
\end{prop}

\begin{proof}
{\bf Claim 1}. The only subgroups isomorphic to $A{\rtimes}{\Z}$ with $A$ finite cyclic are isomorphic to ${\Z}/2{\Z}{\times}{\Z}$.

One can easily check that $\< px, b\>\cong \Z/2{\Z}\times\Z$ and so
$\Z/2{\Z}\times \Z$ is a subgroup group of $G$.

Now the only elements of order three in $\hol(F_2)$ are conjugates of $y^2\t_b^{-1}$ or 
$y^2\t_b^{-1}a$. Assume that there is a subgroup of $G$ isomorphic to $\Z/3\Z \rtimes \Z$. 
Then, conjugating if necessary, we may assume that there is an element of infinite order in $G$, say $g$, such that $g^{-1}(y^2\t_b^{-1}a^s)g=(y^2\t_b^{-1}a^s)^{\pm 1}$ with $s\in \{ 0,1\}$. Based on the decomposition (\ref{simplehol}) we see that the above relation implies that $y^{2}\t_b^{-1}a^s$
stabilizes both vertices $F_2\rtimes C$ and $g^{-1}(F_2\rtimes C)$ and hence the
path that connects them. So it belongs to the edge stabilizer $F_2\rtimes D$, unless $g\in F_2\rtimes C$.
But if $g\in F_2\rtimes D$ we have a contradiction since by decomposition (\ref{furtherhol}),  
$F_2\rtimes D$ cannot contain elements of order $3$. Now if $g\in F_2\rtimes C$,
then based again on the decomposition (\ref{furtherhol}) of $F_2\rtimes C$, we have that 
$g$ stabilizes both $F_2\rtimes D_3$ and $g^{-1}(F_2\rtimes D_3)$ and so it belongs to $F_2\rtimes \Z/2\Z$, a further contradiction, unless again $g\in F_2\rtimes D_3$. Finally, by the decomposition of
$F_2\rtimes D_3= D_3*_{\Z/2\Z}D_3$ we have again that $g$ has to be an element of either of the
two $D_3$ vertices and hence an element of finite order.

We shall show now that $G$ cannot contain subgroups isomorphic to $\Z/4{\Z}\rtimes \Z$. Assume
that $G$ contains such a subgroup, say $A$. Then, $A$ is generated by a conjugate of
$x$, since the conjugacy class of $x$ is the only class of elements of order $4$, and
by an element $g$ of $G$. Using conjugation if necessary, we may assume the element of
order $4$ in $A$ is $x$. Now let $g\in G$ such that $\< x, g\>\cong \Z/4\Z\rtimes\Z.$ 
Then $g^{-1}xg=x^{\pm 1}$. Let $G$ act to the tree that corresponds to the decomposition 
(\ref{simplehol}). Then, due to the above relation, $x$ stabilizes both $F_2\rtimes B$ and 
$g^{-1}(F_2\rtimes B)$ and so it stabilizes the path between the two vertices. 
Hence, $x\in F_2\rtimes D$ a contradiction, unless $g\in F_2\rtimes B$. Moreover, using the decomposition (\ref{furtherhol}) we see that $g$ can only be an element of $F_2\rtimes D_4$ and using the 
fact that $F_2\rtimes D_4=D_4 *_{\Z/2\Z} D_2$ it can only be an element of $D_4$ and so
is of finite order. This completes the proof of claim 1.

\vspace{18pt}\noindent
{\bf Claim 2}.  There are no subgroups of $G$ of the form $D_2 \rtimes \Z$ or of the form 
$D_3\rtimes \Z$.

Up to conjugacy, the possible $D_2$ in $G$ are $\< x^2,p\>$, $\< px, x^2\>$, $\< px, x^2b\>$,
$\< px, x^2\t_b\>$, $\<px, x^2\t_bb^{-1}\>$ and $\< pxa, x^2\t_bb^{-1}\>$. 
Now notice that all last five, $\< px, x^2\>$, $\< px, x^2b\>$,
$\< px, x^2\t_b\>$, $\<px, x^2\t_bb^{-1}\>$ and $\< pxa, x^2\t_bb^{-1}\>$
appear only once in the graph of groups decomposition (\ref{furtherhol}) of $\hol(F_2)$, as vertex groups.
Moreover, in all five, none of the generators is
conjugate to the other, i.e. there are no $g\in G$
such that $gpxg^{-1}=x^2$ or $gpxg^{-1}=x^2b$ or
$gpxg^{-1}=x^2\t_bb^{-1}$ or $gpxg^{-1}=x^2b$ or
$gpxag^{-1}=x^2\t_bb^{-1}$ by Lemma \ref{lem-finite}.
Hence, a relation of the form $gD_2g^{-1}=D_2$ implies (repeating again the argument of claim 1) that $g$ is an element of finite order. So the only possibility for a semidirect
product $D_2\rtimes \Z$ lies with $\< x^2,p\>$.

So assume that there is an element $g\in G$ such that $\< g, p, x^2\>=D_2\rtimes\Z$. Then,
since $p$ and $x^2$ are not conjugates and $px^2$ is conjugate to
$p$, the action of $g$ is either $gpg^{-1}=p$ and $gx^2g^{-1}=x^2$, or $gpg^{-1}=px^2$ and $gx^2g^{-1}=x^2$. 
Let us concentrate to the relation $gpg^{-1}=p$. 
Given the normal form of $g=p^{r}u(x,y)x^{2s}w(\t_a,\t_b)z(a,b)$ we have that
$$p^{r}u(x,y)x^{2s}w(\t_a,\t_b)z(a,b) p z^{-1}(a,b)w^{-1}(\t_a,\t_b)x^{-2s}u^{-1}(x,y)p^{-r}=p.$$
The above relation implies the existence of the following relation in
$\gl$:
$$P^rU(X,Y)X^{2s}PX^{2s}U^{-1}(X,Y)P^r=P$$ which is equivalent to
$$UPU^{-1}=P.$$ By the normal form for the elements of $\gl$, we have that $U$ is of the 
form $U=XY^{e_1}\ldots XY^{e_k}$ with $e_i\in \{\pm 1\}$. So the word $UPU^{-1}$
becomes  
$$XY^{e_1}\ldots XY^{e_k} P Y^{-e_k}X^{-1}\ldots Y^{-e_1}X^{-1}=$$
$$XY^{e_1}\ldots XY^{e_k} Y^{e_k}X\ldots Y^{e_1}X\cdot P=$$
\begin{eqnarray}
X^2\cdot XY^{e_1}\ldots XY^{e_{k-1}}X\cdot Y^{-1}\cdot XY^{e_{k-1}}\ldots Y^{e_1}X\cdot P & {\rm if} & e_k=1\nonumber\\
X^2\cdot XY^{e_1}\ldots XY^{e_{k-1}}X\cdot Y   \cdot XY^{e_{k-1}}\ldots Y^{e_1}X\cdot P & {\rm if} & e_k=-1\nonumber\\
XY^{e_1}\ldots XY^{\pm 1} X\ldots Y^{e_1}XP & {\rm if} & e_k=0\;\;\text{and}\;\;e_{k-1}=\mp 1.\nonumber
\end{eqnarray}
In all cases, the relation $UPU^{-1}=P$ is impossible, since, after deletions of $P$, the remaining word is reduced as written so is never trivial, unless $U=1$. 
Hence $u=1$ and so 
the only possible $g=p^rx^{2s}w(\t_a,\t_b)z(a,b)$. Then the relation $gpg^{-1}=p$ gives
$$x^{2s}w(\t_a,\t_b)z(a,b) p z^{-1}(a,b)w^{-1}(\t_a,\t_b)x^{2s}=p$$ i.e.
$w(\t_b,\t_a)w^{-1}(\t_a,\t_b)=1$. One can easily see that
if $w(\t_a,\t_b)$ is a reduced word in $\inn(F_2)\cong F_2$
then the word $w(\t_b,\t_a)$ is reduced and the 
word $w(\t_b,\t_a)w^{-1}(\t_a,\t_b)$ is reduced and cyclically reduced as written. Hence, a relation 
$w(\t_b,\t_a)w^{-1}(\t_a,\t_b)=1$ is impossible 
unless $w=1$. That implies $z(b,a)z^{-1}(a,b)=1$ and again $z=1$. 
Then $g=p^rx^{2s}$ which has finite order for all possible $r,s$.  

Let us now examine the possibility $gpg^{-1}=px^2$. This implies that 
$$p^{r}u(x,y)x^{2s}w(\t_a,\t_b)z(a,b) p z^{-1}(a,b)w^{-1}(\t_a,\t_b)x^{-2s}u^{-1}(x,y)p^{-r}=px^2.$$
Projection to $\gl$ gives
$$P^rU(X,Y)X^{2s}PX^{2s}U^{-1}(X,Y)P^r=PX^2$$ which is equivalent to
$$UPU^{-1}=PX^2.$$ Performing the same analysis as above for $UPU^{-1}$
we get that the only possibility is $U=X$, hence $u=x$. Then $g=p^rxx^{2s}w(\t_a,\t_b)z(a,b)$
and so the relation $p^rxx^{2s}w(\t_a,\t_b)z(a,b)pz^{-1}(a,b)w^{-1}(\t_a,\t_b)x^{2s}x^{-1}p^r=px^2$ implies 
again that $w(\t_b,\t_a)w^{-1}(\t_a,\t_b)=1$ which possible if and only if $w=1$ and 
so $z(b,a)z^{-1}(a,b)=1$ which is possible if and only if $z=1$. 

Finally, one can easily check that existence of $g\in G$ such that $gD_3g^{-1}=D_3$ can only
occur for $g$ of finite order (using again the previous arboreal argument), so we have that
no subgroup of the form $D_3\rtimes\Z$. 

The only case left is subgroups isomorphic to $D_4\rtimes\Z$. It is easy to see that then it will contain subgroups isomorphic to $D_2\rtimes\Z$, which is impossible. 
\end{proof}

\section{Vanishing of lower $K$-theory of $\hol(F_2)$}

We will prove the main result of the paper. For a group $G$, we write
$$\text{Wh}_q(G) = \left\{
\begin{array}{ll}
\text{Wh}(G), & \text{if}\; q = 1,\\
\tilde{K_0}({\Z}G), & \text{if}\; q = 0, \\
K_q({\Z}G), & \text{if}\; q < 0.
\end{array}\right.$$

\begin{thm}
Let $\Gamma < \hol(F_2)$. Then 
for all $q \le 1$, 
$\text{Wh}_q(\Gamma ) = 0$.
\end{thm}

\begin{proof}

We will show the theorem for $G = \hol(F_2)$. The proof for $\aut(F_2)$ is similar.
By Proposition \ref{prop-hol}, $G$ satisfies the FIC. Let $\Gamma < \hol(F_2)$. Then by \cite{fj1},
$\Gamma$ also satisfies the FIC. 
Thus the maps
$$H_q^G(E_{\fbc}\Gamma ; \K\Z^{-\infty}) \to \text{Wh}_q(\Z \Gamma ), \qquad q \le 1,$$
are isomorphisms. There is a spectral sequence that computes the left hand side of such an isomorphism:
$$E^2_{i,j} = H_i^G(E_{\fbc}\Gamma; \text{Wh}_j(V)) \Longrightarrow \text{Wh}_{i+j}(\Gamma),$$
where $V$ is in $\fbc$. Now, by the decomposition of $\hol(F_2)$ and Proposition \ref{prop-fbc}:
\begin{enumerate}
\item If $V$ is finite, $V$ will isomorphic to one of the following groups: ${\Z}/2{\Z}$, ${\Z}/3{\Z}$, ${\Z}/4{\Z}$, $D_2$, $D_4$. But in this case from the lists in \cite{ao} and \cite{lmo}
$$\text{Wh}_q(V) = 0, \qquad \text{for}\;\; q \le 1.$$
\item If $V$ is infinite, then $V = {\Z}/2{\Z}{\times}{\Z}$. Using the Bass--Heller--Swan Formula and the calculations of the Nil-groups in \cite{bass}, we have that:
$$\text{Wh}_q(V) = 0, \qquad \text{for}\;\; q \le 1.$$
\end{enumerate}
Thus $\text{Wh}_q(\Gamma ) = 0$ for all $q \le 1$.
\end{proof}

\end{document}